\input amstex.tex
\documentstyle{amsppt}
\magnification=1200

\loadbold

\def\Co#1{{\Cal O}_{#1}}
\def\Fi#1{\Phi_{|#1|}}
\def\fei#1{\phi_{#1}}

\def\lrw{\longrightarrow}

\def\Bbbp1{{\Bbb P}^1}
\def\simlin{\sim_{\text{\rm lin}}}

\def\Div{\text{\rm Div}}

\def\dim{\text{\rm dim}}
\def\lrceil#1{\ulcorner {#1}\urcorner}

\TagsOnRight
\pagewidth{14cm}
\pageheight{21cm}

\topmatter
\title
THE RELATIVE PLURICANONICAL STABILITY FOR 3-FOLDS OF GENERAL TYPE
\endtitle
\author   Meng Chen \endauthor
\address \ \ \ \ \ \ \ \ 
\newline Department of Applied Mathematics,\ Tongji University,\ 
 Shanghai,\ 200092,\ P.R. CHINA
\newline E-mail: {\bf chenmb\@online.sh.cn}
\newline \ \ \ \ \ \ 
\newline The Abdus Salam International Centre for Theoretical Physics,
\  Trieste, ITALY
\endaddress
\thanks The author was supported in part by the Abdus Salam International 
Centre for Theoretical Physics and the National Natural Science Foundation of 
China
\endthanks
\abstract
The aim of this paper is to improve a theorem of J\'anos Koll\'ar by a 
different method. For a given smooth
Complex projective threefold $X$ of general type, suppose the plurigenus $P_k(X)\ge 2$,
Koll\'ar proved that the $(11k+5)$-canonical map is birational.
Here we show that either the $(7k+3)$-canonical map or the $(7k+5)$-canonical 
map is birational and the $(13k+6)$-canonical
map is stably birational onto its image. If  $P_k(X)\ge 3$, then
the $m$-canonical map is  birational for $m\ge 10k+8$.
In particular, $\fei{12}$ is birational when $p_g(X)\ge 2$ and 
$\fei{11}$ is birational when $p_g(X)\ge 3$.
\endabstract
\endtopmatter

\document
\baselineskip 12.5pt
\leftheadtext{Meng Chen}
\rightheadtext{The Relative Pluricanonical Stability}

\head {Introduction}\endhead
Let $X$ be a smooth projective 3-fold of general type defined over ${\Bbb C}$
and  denote
by $\phi_m$ the $m$-canonical map of $X$, which is the rational map associated
with the linear system $|mK_X|$. 
Let $P_k(X):=h^0(X,\Co{X}(kK_X))$ for any positive integer $k$, we usually call
$P_k(X)$ {\it the k-th plurigenus} of $X$ which is a birational invariant.
For a given positive integer $m_0$, we say that $\fei{m_0}$ is {\it stably 
birational} if $\fei{m}$ is birational onto its image for all $m\ge m_0$.
Since the Kodaira dimension $\text{kod}(X)=3$, $\fei{m}$ is 
birational for $m\gg 0$.
In this paper, we consider the following
\proclaim{Problem} Suppose $P_k(X)\ge 2$,
for which value $m_0(k)$, does $|m_0(k)K_X|$ define a stably birational map onto
its image?
\endproclaim
In 1986, Koll\'ar (\cite{5, Corollary 4.8}) first gave an effective result 
and proved that the $(11k+5)$-canonical map is birational if $P_k(X)\ge 2$.
However, his method can not tell whether $\fei{m}$ is still birational for all
$m>11k+5$. On the other hand, it seems to us that the number $11k+5$ is not 
the optimal one. This paper aims to present a better result as the following

\proclaim{Main Theorem} Let $X$ be a nonsingular projective threefold of general type and
 suppose $P_k(X)\ge 2$, then

(i) either $\fei{7k+3}$ or $\fei{7k+5}$ is birational onto its image;

(ii) $\fei{13k+6}$ is stably birational onto its image;

(iii) $\fei{10k+8}$ is stably birational providing that $P_k(X)\ge 3$.

\noindent In particular, if $p_g(X)\ge 2$, then $\fei{m}$ is birational for all
$m\ge 12$; if $p_g(X)\ge 3$, then $\fei{m}$ is birational for all $m\ge 11$.
\endproclaim

Noting that the main obstacle which prevents Koll\'ar's method from getting 
a better bound is the case when $X$ admits a rational pencil of certain surfaces
of general type, we mainly take a special study to this 
situation in an alternative way. First we build some birationality criteria
for adjoint systems on a surface of general type, then we reduce the problem 
to surface case while finding suitable divisors on the threefold whose 
restrictions to the surface satisfy those criteria. The Kawamata-Viehweg 
vanishing theorem plays a key role throughout our argument.

\definition{Definition} Let $X$ be a normal projective variety and $D$ be a Weil
divisor on $X$. Denote by $\Fi{D}$ the natural rational map defined by the linear 
system $|D|$. $|D|$ is called {\it base point free} if it has neither fixed 
components nor base points. 

If $|L|$ is a linear system on $X$ without fixed components and $h^0(X, L)\ge 2$, 
we mean {\it a general irreducible element} $S$ of $|L|$ as follows:

(1) if $\dim\Fi{L}(X)\ge 2$, then $S$ is a general member of $|L|$. 

(2) if $\dim\Fi{L}(X)=1$, then $L$ is linearly equivalent to a  union 
of distinct reduced irreducible divisors of the same type. 
Explicitly, $L\simlin \sum S_i$. We mean $S$ a general $S_i$.

$X$ is called {\it minimal} if the canonical divisor $K_X$ is nef, i.e. 
$K_X\cdot C\ge 0$ for all proper curve $C\subset X$. 

$X$ is said to {\it be of general type} if the Kodaira dimension 
$\text{kod}(X)=\dim(X)$.

$X$ is said to {\it have only terminal singularities} according to Reid (\cite{7})
if the following two conditions hold:

(i) for some integer $r\ge 1$, $rK_X$ is Cartier;

(ii) for some resolution $f:Y\lrw X$, $K_Y=f^*(K_X)+\sum a_iE_i$ for 
$0<a_i\in{\Bbb Q}$ for all $i$, where the $E_i$ vary all the exceptional divisors 
on $Y$. 
\enddefinition

\remark{\smc Acknowledgment}
I would like to thank the Abdus Salam International Centre for Theoretical
Physics, Trieste, Italy for support during my visit. I also wish to thank
Prof. L. Ein and Prof. M. S. Narasimhan for encouragement. Thanks are also due to
Prof. J. Koll\'ar and Prof. M. Reid, who privately gave me much more useful
information on the problem considered in this paper. During the preparation of this 
note, I took visits to both Universit\`a di Roma II and Universit\`a 
degli Studi di Padova to the invitation of Prof. Ciro Ciliberto and Prof. Ezio
Stagnaro whom I take this chance to thank for their hospitality.
Finally, the referee's skillful comments greatly helped me to organize this paper
into the present form.
\endremark

\head {1. Preparation}\endhead
Throughout our argument, the Kawamata-Viehweg vanishing theorem is always employed
as a much more effective tool. We use it in the following form.
\proclaim{Vanishing Theorem} (\cite{3} or \cite{10}) Let $X$ be a nonsingular
complete variety, $D\in\Div(X)\otimes{\Bbb Q}$. Assume the following two conditions:

(1) $D$ is nef and big;

(2) the fractional part of $D$ has supports with only normal crossings.

\noindent Then $H^i(X,\Co{X}(\lrceil{D}+K_X))=0$ for $i>0$, where $\lrceil{D}$
is the round-up of $D$, i.e. the minimum integral divisor with $\lrceil{D}-D\ge 0$.
\endproclaim

Another important principle that is tacitly used throughout the text is due to Tankeev 
(\cite{9}). Explicitly, on a smooth projective variety $X$, if we have a base point
free system $|M|$ and an effective divisor $D$, we want to study the birationality of the 
map $\Fi{D+M}$. Now let $S$ be a general irreducible element of $|M|$, then $S$ is 
a smooth divisor on $X$ by Bertini's theorem. Suppose we have known that 
$\Fi{D+M}$ can distinguish general irreducible elements and that 
$\Fi{D+M}\bigm|_S$ is birational, then Tankeev's principle implies the birationality
of $\Fi{D+M}$.
 
\proclaim{Lemma 1.1} (\cite{8, Corollary 2}) Let $S$ be a nonsingular algebraic surface,
$L$ be a nef divisor on $S$, $L^2\ge 10$ and let $\phi$ be a map defined by
$|L+K_S|$. If $\phi$ is not birational, then $S$ contains a base point free
pencil $E'$ with $L\cdot E'=1$ or $L\cdot E'=2$.
\endproclaim

\proclaim{Lemma 1.2} Let $S$ be a nonsingular projective surface of general type,
suppose $L$ be a divisor with $h^0(S,L)\ge 2$, then $h^0(S, K_S+L)\ge 2.$
In particular, if $\chi(\Co{S})\ge 3$, then $h^0(S, K_S+L)\ge 4$.
\endproclaim
\demo{Proof} 
Taking a general irreducible element $C$ in the moving part of $|L|$, then $C$ is a
nef divisor, $C\le L$ and $C$ is a curve of genus $\ge 2$. By R-R on the surface 
$S$, we have
$$h^0(S,K_S+L)\ge h^0(S,K_S+C)\ge\frac{1}{2}(K_S\cdot C+C^2)+\chi(\Co{S}).$$
It is easy to get the result.
\qed
\enddemo

\proclaim{Lemma 1.3} Let $S$ be a nonsingular projective surface of general 
type, $L$ be a nef divisor, $L^2\ge 3$ and $\dim\Fi{L}(S)=2$, then 
$|K_S+2L|$ gives a birational map.
\endproclaim
\demo{Proof} We have $(2L)^2\ge 12$. If $\Fi{K_S+2L}$ is not birational, then 
according to Lemma 1.1, there is a base point free pencil $E'$ such that 
$2L\cdot E'\le 2$, i.e. $L\cdot E'=1$. Since $\dim\Fi{L}(S)=2$ and $E'$ is 
a curve of genus $\ge 2$, we see that $L\cdot E'\ge 2$, a contradiction.
\qed
\enddemo

\proclaim{Lemma 1.4} Let $S$ be a nonsingular projective surface of general type,
$L_i$ be a divisor on $S$ such that $\dim\Fi{L_i}(S)\ge i$ for $i=1,\ 2$, then
$|K_S+2L_2+L_1|$ gives a birational map.
\endproclaim
\demo{Proof} Modulo blowing-ups,
 we can suppose that the $|L_i|$ be  base point free  for $i=1,\ 2$. This 
 means that $L_2$ is nef and big and that $L_1$ is nef.

If the system $|L_2|$ gives a birational map, then so does $|K_S+2L_2+L_1|$,
because $K_S+L_1$ is effective by Lemma 1.2.

Otherwise, we have $L_2^2\ge 2$. Now we have $(2L_2+L_1)^2\ge 12$. If
$|K_S+2L_2+L_1|$ does not give a birational map, then, by Lemma 1.1, there is
a free pencil $E'$ on $S$ such that
$$(2L_2+L_1)\cdot E'\le 2.$$
This means $L_2\cdot E'=1$. Note that $E'$ is a curve of genus $\ge 2$ and
$|L_2|$ gives a generically finite map. The Riemann-Roch theorem on
the curve $E'$ tells that
$\deg(L_2|_{E'})\ge 2$. We have derived a contradiction.
\qed\enddemo

\proclaim{Lemma 1.5} Let $X$ be a nonsingular projective 3-fold of general
type. Suppose  $L_i$ be a divisor on $X$ such that $\dim\Fi{L_i}(X)\ge i$ for
$i=1,\ 2,\ 3$, then $|K_X+2L_3+L_2+L_1|$ gives a birational map.
\endproclaim
\demo{Proof}
Take a birational modification $\pi:X'\lrw X$ according to Hironaka such that
the $|\pi^*(L_i)|$ are all base point free for $i>0$. On $X'$, we can study the
system $|K_{X'}+2\pi^*(L_3)+\pi^*(L_2)+\pi^*(L_1)|$. Let $M_i$ be the moving
part of $|\pi^*(L_i)|$, we have
$$|K_{X'}+2M_3+M_2+M_1|\subset |K_{X'}+2\pi^*(L_3)+\pi^*(L_2)+\pi^*(L_1)|.$$
Therefore, for simplicity, we can suppose from the beginning that the $|L_i|$ are
base point free  on $X$. So $L_3$ is nef and big under this assumption.

Step 1. Verifying that $K_X+2L_3+L_2$ is effective.

We have $\dim\Fi{L_2}(X)\ge 2$. So a general member $S\in|L_2|$ is a nonsingular
projective surface of general type. Using the vanishing theorem to the exact sequence
$$0\lrw\Co{X}(K_X+2L_3)\lrw\Co{X}(K_X+2L_3+S)\lrw\Co{S}(K_S+2L_3|_S)\lrw 0,$$
we get the surjective map
$$H^0(X, K_X+2L_3+S)\lrw H^0(S, K_S+2L_3|_S)\lrw 0.$$
{}From Lemma 1.2, we know $K_S+2L_3|_S$ is effective, so is $K_X+2L_3+L_2$.

Step 2. Reduction to surface case.

Taking a 1-dimensional sub-system  of $|L_1|$, then this system defines a rational
map onto ${\Bbb P}^1$. Taking further blowing-up if necessary, we can also suppose
that this system defines a morphism $f:X\lrw {\Bbb P}^1$. 
Taking the Stein factorization of $f$, one obtains a derived fibration 
$g:X\lrw C$. A general fibre of $f$ can be written as a disjoint union $\sum F_i$.
Let $F$ be a general fibre of $g$, then it is a nonsingular projective surface
of general type and we have $F\le L_1$. Now considering the system
$|K_X+2L_3+L_2+\sum F_i|$, it can distinguish general fibres of $g$
because of  $K_X+2L_3+L_2$ is effective and $2L_3+L_2$ is nef and big.
 Using the vanishing theorem again, we have
$$|K_X+2L_3+L_2+\sum F_i|\bigm|_F=|K_F+2L_3'+L_2'|,$$
where $L_3':=L_3|_F$ and $L_2':=L_2|_F$. Lemma 1.4 tells that the right system gives
a birational map, so does $|K_X+2L_3+L_2+L_1|$. The proof is completed.
\qed\enddemo

\proclaim{Lemma 1.6} Let $X$ be a nonsingular variety of dimension $n$,
$D\in\Div (X)\otimes{\Bbb Q}$ be a ${\Bbb Q}$-divisor on $X$. Then we have
the following:

(i) if $S$ is a smooth irreducible divisor on $X$, then
$\lrceil{D}|_S\ge \lrceil{D|_S}$;

(ii) if $\pi: X'\lrw X$ is a birational morphism, then
$\pi^*(\lrceil{D})\ge \lrceil{\pi^*(D)}$.
\endproclaim
\demo{Proof}
We can write $D$ as $G+\sum_{i=1}^t a_iE_i$, where $G$ is a divisor,
the $E_i$ are effective divisors for each $i$ and $0<a_i<1$, $\forall\ i$.
So we only have to prove the lemma for effective ${\Bbb Q}$-divisors.
That is easy to check.
\qed\enddemo

\proclaim{Lemma 1.7} 
Let $X$ be a nonsingular projective threefold of general type. Let $D$ be a divisor
on $X$ with $h^0(X, D)\ge 2$ and suppose $|D|$ has no fixed components. Denote
by $F$ a general irreducible element of $|D|$. If $L$ is another divisor such that
$\dim\Fi{L}(F)\ge 1$, then $mK_X+L+D$ is effective and $\dim\Fi{mK_X+L+D}(F)\ge 1$
for all $m\ge 2.$ 
\endproclaim
\demo{Proof} 
According to the 3-dimensional MMP (\cite{4} and \cite{6}), $X$ has a minimal
model $X_0$ which is normal projective with only ${\Bbb Q}$-factorial terminal
singularities. Let $\alpha:X\dashrightarrow X_0$ be the contraction which 
is a rational map. Take a common resolution $X'$ with $\pi':X'\lrw X$ and
$\pi:X'\lrw X_0$ such that $\pi=\alpha\circ\pi'$ and that

(1) both $|{\pi'}^*(L)|$ and $|{\pi'}^*(D)|$ have no base points (they may have 
fixed components);

(2) $\pi^*(K_{X_0})$ has supports with only normal crossings.

This is possible because of Hironaka's big theorem. Since
${\pi'}^*(mK_X+L+D)\le mK_{X'}+{\pi'}^*(L)+{\pi'}^*(D)$ and
$${\pi'}_*\Co{X'}(mK_{X'}+{\pi'}^*(L)+{\pi'}^*(D))=
\Co{X}(mK_X+L+D)={\pi'}_*{\pi'}^*\Co{X}(mK_X+L+D),$$
then $h^0\bigl(X',{\pi'}^*(mK_X+L+D)\bigr)=h^0\bigl(X',mK_{X'}+{\pi'}^*(L)+{\pi'}^*(D)\bigr),$
so
$$\Fi{{\pi'}^*(mK_X+L+D)}\ \text{and}\ \Fi{mK_{X'}+{\pi'}^*(L)+{\pi'}^*(D)}$$
have the same behavior. Let $S$ be a general irreducible element of the moving part of 
$|{\pi'}^*(D)|$, then $\dim\Fi{{\pi'}^*(L)}(S)\ge 1$ by assumption. Therefore it is 
sufficient to show 
$$\dim\Fi{mK_{X'}+{\pi'}^*(L)+{\pi'}^*(D)}(S)\ge 1$$
for $m\ge 2$. Let $H$ be the moving part of $|{\pi'}^*(L)|$, then $H$ is nef since
$|H|$ is base point free.
 We have
$$|K_{X'}+\lrceil{(m-1)\pi^*K_{X_0}}+H+S|\subset|mK_{X'}+{\pi'}^*(L)+{\pi'}^*(D)|.$$
The Kawamata-Viehweg vanishing theorem gives
$$\align
&|K_{X'}+\lrceil{(m-1)\pi^*K_{X_0}}+H+S|\bigm|_S\\
&=\bigm|K_S+\lrceil{(m-1)\pi^*K_{X_0}}|_S+M\bigm|\supset|K_S+\lrceil{B}+M|,
\endalign$$
where $B:=(m-1)\pi^*K_{X_0}|_S$ is nef and big on $S$ and  $M:=H|_S$.
From the assumption, we have $h^0(S, M)\ge 2$. Choosing a 1-dimensional sub-system
$|C|$ in $|M|$, modulo blowing-ups, we can suppose $|C|$ be base point free.
Also from the vanishing theorem, we have
$$|K_S+\lrceil{B}+C|\bigm|_C=|K_C+D|,$$
where $D:=\lrceil{B}|_C$ is a divisor on the curve $C$ with positive degree
since $D\ge \lrceil{B|_C}$ by Lemma 1.6(i).
Because $g(C)\ge 2$, we have $h^0(K_C+D)\ge 2$. This means $|K_C+D|$ gives a
generically finite map and
$$\dim \Fi{K_S+\lrceil{B}+C}(C)=1$$
 thus $K_{X'}+\lrceil{(m-1)\pi^*K_{X_0}}+{\pi'}^*(L)+{\pi'}^*(D)$
is effective and the image of $S$ through the map defined by this divisor is
at least $1$. The proof is completed.
\qed
\enddemo

\head {2. Proof of the main theorem} \endhead

\subhead 2.1 Basic formula\endsubhead
Let $X$ be a nonsingular projective threefold, $f: X\lrw C$ be a fibration onto
a nonsingular curve $C$. {}From the spectral sequence:
$$E_2^{p,q}:=H^p(C,R^qf_*\omega_X)\Rightarrow E^n:=H^n(X,\omega_X),$$
we get by direct calculation that
$$h^2(X, \Co{X})=h^1(C,f_*\omega_X)+h^0(C,R^1f_*\omega_X),$$
$$q(X):=h^1(X,\Co{X})=b+h^1(C,R^1f_*\omega_X),$$
where $b$ denotes the genus of $C$.

\subhead 2.2 Review of Koll\'ar's technique\endsubhead
Let $X$ be a smooth projective 3-fold of general type and suppose $P_k(X)\ge 2$.
Choose a $1$-dimensional sub-system of $|kK_X|$ and replace $X$ by a
birational model $X'$ where
this pencil defines a morphism $g:X'\lrw {\Bbb P}^1$. (For simplicity, we can
suppose $X'=X$).  Let $S$ be a general irreducible element of this pencil, 
then a general fibre of $g$ is a disjoint union of some surfaces with the same 
type as $S$ and  
 $S$ is a smooth projective surface of general type. Let
$t=k(2p+1)+p$. Then $H^0(\omega_X^t)=
H^0({\Bbb P}^1, g_*\omega_X^t)$ and we have an injection
$\Cal{O}(1)\hookrightarrow g_*\omega_X^k$, and hence an injection
$\Cal{O}(2p+1)\hookrightarrow g_*\omega_X^{k(2p+1)}$. This gives an injection
$$\Cal{O}(2p+1)\otimes g_*\omega_X^p\hookrightarrow g_*\omega_X^t,$$
where $\Cal{O}(2p+1)\otimes g_*\omega_X^p=
\Cal{O}(1)\otimes g_*\omega_{X/{\Bbb P}^1}^p$.
Now it is well-known that $g_*\omega_{X/{\Bbb P}^1}^p$ is a sum of line bundles
of non-negative degree on ${\Bbb P}^1$. If $p\ge 5$, the local sections of
$g_*\omega_X^p$ give a birational map for $S$, and all these extend to global
sections of $\Cal{O}(2p+1)\otimes g_*\omega_X^p$. Moreover its sections separate
the fibres from each other, hence $\phi_t$ is a birational map for $X$.

{}From the above method, according to \cite{1} and \cite{11}, we have

(1) $\phi_{5k+2}$ is generically finite for $X$ if $S$ is not a surface with
$p_g(S)=q(S)=0$ and $K_{S_0}^2=1$, where $S_0$ is the minimal model of $S$.
Otherwise, we have at least $\dim\phi_{5k+2}(X)\ge 2$;

(2) $\phi_{7k+3}$ is birational for $X$ if $S$ is not a surface with
$$(K_{S_0}^2, p_g(S))=(1,2)\ \text{or}\ (2,3).$$

\subhead 2.3 Proof of the main theorem\endsubhead
According to the 3-dimensional MMP, we can suppose $X$ be a minimal model
with at worst ${\Bbb Q}$-factorial terminal singularities. This means that
$K_X$ is a nef and big ${\Bbb Q}$-divisor. We begin from a minimal model
in order to make use of the Kawamata-Viehweg vanishing theorem. 

\proclaim{Theorem 2.3.1} Let $X$ be a nonsingular projective 3-fold of general type and
suppose $P_k(X)\ge 2$, then either $\fei{7k+3}$ or $\fei{7k+5}$ is birational.
\endproclaim
\demo{Proof}
Suppose $X$ be a minimal model with at worst ${\Bbb Q}$-factorial terminal 
singularities. Choose a 1-dimensional sub-system $\Lambda$ of $|kK_X|$ and take 
a birational modification $\pi:X'\lrw X$ such that 

(i) $X'$ is nonsingular;

(ii) $\pi^*\Lambda$ gives a morphism;

(iii) the fractional part of $\pi^*(K_X)$ has supports with only normal 
crossings.

This is possible because of Hironaka's big theorem. Set 
$g_1:=\Phi_{\Lambda}\circ\pi$ 
and let
$X'\overset f_1\to\lrw W_1\overset s_1\to\lrw \Bbbp1$ be the Stein factorization
of $g_1$. Denote $b:=g(W_1)$, the geometric genus of the curve $W_1$.

If $b>0$, then the moving part of $\Lambda$ is base point free. Let $\sum S_i$ 
be the moving part of $\Lambda$, then $\sum S_i\le kK_X$ and a general $S_i$ is a 
smooth projective surface of general type, since the singularities on $X$ are 
isolated. Using Kawamata's vanishing theorem 
(\cite{4}) to ${\Bbb Q}$-Cartier Weil divisors on minimal threefold $X$, 
we see that $|(a+1)K_X+\sum S_i|$ can distinguish general $S_i$ for $a>0$ and
$$H^0(X, (a+1)K_X+\sum S_i)\lrw \oplus H^0(S_i, (a+1)K_{S_i})$$
is surjective. Therefore it is obvious that $\fei{m}$ is effective whenever 
$m\ge k+2$, generically finite whenever $m\ge 2k+2$, birational whenever 
$m\ge 2k+4$.

So, from now on, we can suppose that $b=0$. We have a fibration 
$f_1:X'\lrw \Bbbp1$. Let $F$ be a general fibre of $f_1$. By virtue of 2.2(2),
we can suppose that $F$ is a surface with invariants $(K_{F_0}^2,p_g(F))=(1,2)$ or
$(2,3)$, where $F_0$ is the minimal model of $F$. $F$ is the moving part 
of $\pi^*\Lambda$ and $F\le_{\Bbb Q} \pi^*(kK_X)$. 
We automatically have $q(F)=0$. First we study the system
 $|K_{X'}+\lrceil{k\pi^*(K_X)}+F|.$
 For a general fibre $F$, the vanishing theorem gives that
$$|K_{X'}+\lrceil{k\pi^*(K_X)}+F|\bigm|_F=\bigm|K_F+\lrceil{k\pi^*(K_X)}|_F\bigm|,$$
where $\lrceil{k\pi^*(K_X)}|_F$ is effective. This means that $(2k+1)K_{X'}$ is
effective and $\dim\phi_{2k+1}(F)\ge 1$. By Lemma 1.7, we see that
 $mK_{X'}$ is effective and $\dim\phi_m(F)\ge 1$ for $m\ge 3k+3$.

Actually, we have $\dim\phi_{3k+2}(F)=2$. In fact, we have
$$|K_{X'}+\lrceil{(2k+1)\pi^*(K_X)}+F|\bigm|_F\supset\bigm|K_F+M_{2k+1}|_F\bigm|,$$
where $M_{2k+1}$ is the moving part of $|\lrceil{(2k+1)\pi^*K_X}|$. It is easy to check
that $\bigm|K_F+M_{2k+1}|_F\bigm|$ gives a generically finite map because $q(F)=0$ and
$p_g(F)>0$. Thus 
$$\dim\Fi{K_{X'}+\lrceil{(2k+1)\pi^*(K_X)}+F}(F)\ge 2.$$
We have $|K_{X'}+\lrceil{2(3k+2)\pi^*(K_X)}+F|\subset |(7k+5)K_{X'}|.$
$K_{X'}+\lrceil{2(3k+2)\pi^*(K_X)}$ is effective by the above argument. So 
$|K_{X'}+\lrceil{2(3k+2)\pi^*(K_X)}+F|$ can distinguish general fibre $F$. On the
 other hand, the Kawamata-Viehweg vanishing theorem gives 
 $$\align
 |K_{X'}+\lrceil{2(3k+2)\pi^*(K_X)}+F|\bigm|_F&=\bigm|K_F+\lrceil{2(3k+2)\pi^*(K_X)}|_F\bigm|\\
 &\supset |K_F+2L_{3k+2}|,
 \endalign$$
where $L_{3k+2}:=M_{3k+2}|_F$. It is 
sufficient to show that $|K_F+2L_{3k+2}|$ gives a birational map for $F$.
We have already known that $|L_{3k+2}|$ gives a generically finite map for $F$.
Excluding the fixed components of $|L_{3k+2}|$, we can suppose that 
$|L_{3k+2}|$ are moving on the surface $F$. So $L_{3k+2}$ is nef.
If $|L_{3k+2}|$ gives a birational map, then so does $|K_F+2L_{3k+2}|$. Otherwise,
$$L_{3k+2}^2\ge 2(h^0(F,L_{3k+2})-2).$$
Considering the following three natural maps
$$\align
&H^0(X', M_{3k+2})\overset\alpha\to\lrw H^0(F, L_{3k+2})\\
&H^0(X', K_{X'}+\lrceil{(2k+1)\pi^*(K_X)}+F)\overset \beta\to\lrw 
H^0(F,K_F+\lrceil{(2k+1)\pi^*(K_X)}|_F)\lrw 0\\
&H^0\bigl(X', (3k+2)K_{X'}\bigr)\overset \gamma\to\lrw H^0\bigl(F, (3k+2)K_F\bigr)
\endalign$$
where $\beta$ is surjective by the Kawamata-Viehweg vanishing theorem. We see that 
$$\dim_{\Bbb C}\bigl(\text{im}(\alpha)\bigr)=\dim_{\Bbb C}\bigl(\text{im}(\gamma)\bigr)
\ge\dim_{\Bbb C}\bigl(\text{im}(\beta)\bigr)=h^0(F, K_F+D_{2k+1})$$
where $D_{2k+1}:=\lrceil{(2k+1)\pi^*(K_X)}|_F$ and $h^0(F,D_{2k+1})\ge 2$.
So $h^0(F, K_F+D_{2k+1})\ge 4$, according to Lemma 1.2, 
because we have $\chi(\Co{F})\ge 3$ in this case. Thus
$$L_{3k+2}^2\ge 2\bigl(h^0(F,L_{3k+2})-2\bigr)\ge 2\Bigl(\dim_{\Bbb C}\bigl(\text{im}(\alpha)\bigr)-2\Bigr)
\ge 4$$
and then $|K_F+2L_{3k+2}|$ gives a birational map by Lemma 1.3.
So $\fei{7k+5}$ is birational. 

Finally, for all $m\ge 10k+7$, set $t:=m-7k-5\ge 3k+2$, then 
$\dim\fei{t}(F)\ge 1$. 
In particular, $tK_{X'}$ is effective.  So $\fei{m}$ is birational  for all 
$m\ge 10k+7$ in this case.
\qed
\enddemo

\proclaim{Corollary 2.3.1} Let $X$ be an irregular nonsingular 3-fold of general type, suppose
$P_k(X)\ge 2$, then $\phi_{7k+3}$ is birational. Therefore at least $\phi_{143}$
is birational according to Koll\'ar and Fletcher.
\endproclaim
\demo{Proof} In the proof of the last theorem, if $b>0$, then $\fei{m}$ is 
birational for $m\ge 2k+4$. If $b=0$, we can 
use the formula of $q(X)$ to the fibration $f_1:X'\lrw \Bbbp1$.
When $q(X)>0$, then we must have $q(F)>0$. Then $\Fi{3K_F}$ is birational for the 
fibre $F$, so is $\Fi{(7k+3)K_X}$ by 2.2(2). Moreover, we have $P_{20}(X)\ge 2$ for any irregular
3-fold of general type according to Koll\'ar (\cite{5}) and Fletcher (\cite{2}). 
Thus $\phi_{143}$ is birational.
\qed
\enddemo

\proclaim{Theorem 2.3.2} Let $X$ be a nonsingular projective threefold of general type and
suppose $P_k(X)\ge 2$, then $\fei{m}$ is birational for $m\ge 13k+6$.
\endproclaim
\demo{Proof}
Suppose $X$ be a minimal model with at worst ${\Bbb Q}$-factorial terminal 
singularities.
Make a birational modification
$\pi:X'\lrw X$ such that:

(i) $X'$ is nonsingular;

(ii) $|kK_{X'}|$ gives a morphism;

(iii) the fractional part of $\pi^*(K_X)$ has supports with only normal crossings.

Set $g:=\Fi{kK_X}\circ\pi$ and $W':=\overline{\Fi{kK_X}(X)}$. Let
$X'\overset f\to\lrw W\overset s\to\lrw W'$
be the Stein factorization of $g$.

We would like to formulate our proof through two steps as follows.

\noindent{\bf Case 1}. $\dim\phi_k(X)\ge 2$.

Set $kK_{X'}\simlin M_k+Z_k$, where $M_k$ is the moving part and $Z_k$ is the
fixed part. Then a general member $S\in |M_k|$ is an irreducible nonsingular
projective surface of general type.
Write $K_{X'}=\pi^*(K_X)+\sum a_iE_i$, where the $E_i$ are exceptional divisors
for $\pi$, $0<a_i\in {\Bbb Q}$ for each $i$. Obviously, $\lrceil{\pi^*(K_X)}\le
K_{X'}$. Because $h^0(X',\lrceil{\pi^*(kK_X)})=h^0(X', kK_{X'})$, we can see 
that $M_k$ is actually also the moving part of $|\lrceil{\pi^*(kK_X)}|$.
Thus we have
$$\pi^*(kK_X)\ge_{{\Bbb Q}} M_k+\sum b_iE_i,$$
where  $0\le b_i\in {\Bbb Q}$ for each $i$.

We claim that $mK_{X'}$ is always effective for $m\ge 2k+1$. In fact, for any
$t\in {\Bbb Z}^+$, we consider the system $$|K_{X'}+\lrceil{\pi^*((t+k)K_X)}+S|.$$
It is a sub-system of $|(2k+t+1)K_{X'}|$. By the Kawamata-Viehweg vanishing theorem, we have
a surjective map
$$H^0(X', K_{X'}+\lrceil{\pi^*((t+k)K_X)}+S)\lrw
H^0(S, K_S+\lrceil{\pi^*((t+k)K_X)}|_S)\lrw 0.$$
Noting that $\lrceil{\pi^*((t+k)K_X)}\ge \lrceil{\pi^*(tK_X)}+M_k$, also by
Lemma 1.6(i), it is sufficient to show that
$K_S+\lrceil{\pi^*(tK_X)|_S}+M_k|_S$
is effective. 
When $t=0$, then $h^0(S, K_S+M_k|_S)\ge 2$  by Lemma 1.2,
because $h^0(S, M_k|_S)\ge 2$. When $t>0$, choose a
1-dimensional sub-system $|C|$ in the moving part of $\bigm|M_k|_S\bigm|$.  
Modulo blowing-ups, we can suppose $|C|$ be free from base points and 
then $C$ is nef
and $C\le M_k|_S$.
 We have $g(C)\ge 2$. Because
$\pi^*(tK_X)|_S$ is a nef and big ${\Bbb Q}$-divisor on $S$, by the 
 Kawamata-Viehweg vanishing
theorem, we also get a surjective map
$$H^0(S, K_S+\lrceil{\pi^*(tK_X)|_S}+C)\lrw H^0(C, K_C+D)\lrw 0,$$
where $D:=\lrceil{\pi^*(tK_X)|_S}\bigm|_C$ is a divisor on $C$ with positive degree.
Thus $h^0(C, K_C+D)\ge 2$. This leads to the effectiveness of $(2k+t+1)K_{X'}$.
Moreover, actually we have proved that
$\dim\phi_m(S)\ge 1$ for $m\ge 2k+1$.

Now we prove that $\phi_{3k+1}$ is generically finite. Considering the system
$$|K_{X'}+\lrceil{2k\pi^*(K_X)}+M_k|,$$ 
as we have shown in above that $(2k+1)K_{X'}$
is effective, so $|K_{X'}+\lrceil{2k\pi^*(K_X)}+M_k|$ can distinguish general $S$.
 By the Kawamata-Viehweg vanishing theorem, we have
$$|K_{X'}+\lrceil{2k\pi^*(K_X)}+S|\bigm|_S=\bigm|K_S+\lrceil{2k\pi^*(K_X)}|_S\bigm|.$$
We have
$$\bigm|K_S+\lrceil{2k\pi^*(K_X)}|_S\bigm|\supset \bigm|K_S+\lrceil{k\pi^*(K_X)|_S}+M_k|_S\bigm|.$$
Noting that $h^0(S, M_k|_S)\ge 2$,  $K_S+\lrceil{k\pi^*(K_X)|_S}\ge
K_S+M_k|_S$, which is also effective by Lemma 1.2, and $k\pi^*(K_X)|_S$ is a
nef and big ${\Bbb Q}$-divisor on $S$, it is easy to verify that
$\bigm|K_S+\lrceil{k\pi^*(K_X)|_S}+M_k|_S\bigm|$ gives a generically finite map. In fact,
choose a 1-dimensional sub-system $|C|$ in the moving part of $\bigm|M_k|_S\bigm|$.
With the same reason, we can suppose $|C|$ be free from base points.
$\bigm|K_S+\lrceil{k\pi^*(K_X)|_S}+C\bigm|$ can distinguish general $C$, and we have
$$|K_S+\lrceil{k\pi^*(K_X)|_S}+C|\bigm|_C=|K_C+D|,$$
where $D$ is a divisor on $C$ with positive degree. Because $g(C)\ge 2$,
thus $h^0(K_C+D)\ge 2$ and $|K_C+D|$ gives a generically finite map.

Finally, we want to show that $\phi_m$ is birational for  $m\ge 9k+4$.
Let $t:=m-7k-3$, then $t\ge 2k+1$. Denote by $M_{3k+1}$ the moving part of
$|(3k+1)K_{X'}|$ and by $M_t$ the moving part of $|tK_{X'}|$.
We have
$$|K_{X'}+\lrceil{(t+6k+2)\pi^*(K_X)}+M_k|\subset|mK_{X'}|.$$
Because $t+6k+3>2k+1$, $K_{X'}+\lrceil{(t+6k+2)\pi^*(K_X)}$ is effective,
thus the left system in above can distinguish general $S$.
Furthermore, the vanishing theorem gives
$$|K_{X'}+\lrceil{(t+6k+2)\pi^*(K_X)}+M_k|\bigm|_S=|K_S+L|,$$
where $L:=\lrceil{(t+6k+2)\pi^*(K_X)}\bigm|_S\ge 2M_{3k+1}|_S+M_t|_S$. By Lemma 1.4,
$|K_S+L|$ gives a birational map, so does $|mK_{X'}|$.
\smallskip

\noindent{\bf Case 2}. $\dim\phi_k(X)=1$.

In this case, $W$ is a nonsingular curve of genus $b$. Let $F$ be a general
fibre of $f$, then $F$ is an irreducible smooth projective surface of general
 type. We have $M_k\simlin \sum F_i$, where the $F_i$ are fibres of $f$ for
 each $i$.

By a parallel argument as in the proof of Theorem 2.3.1, we see that $\fei{m}$
is birational for $m\ge 2k+4$ if $b>0$. And if $b=0$ while $F$ is a surface with
the invariants $\bigl(K_{F_0}^2, p_g(F)\bigr)=(1,2)$ or $(2,3)$, then $\fei{m}$ is 
birational for $m\ge 10k+7$. 

Otherwise, we use Koll\'ar's method.
 {}From 2.2, we know that $\phi_{7k+3}$ is birational and
$\dim\phi_{5k+2}(X)\ge 2$. Thus, by Lemma 1.7, $mK_{X'}$ is effective
for $m\ge 6k+4$. 
Since we have $|K_{X'}+\lrceil{(5k+2)\pi^*(K_X)}+F|\bigm|_F=|K_F+D|$
where $D:=\lrceil{(5k+2)\pi^*(K_X)}\bigm|_F$ is effective and $h^0(F,D)\ge 2$, we see 
that $K_F+D$ is effective and thus 
$(6k+3)K_{X'}$ is effective. So $\phi_m$ is birational for 
$m\ge 13k+6$, which means that $\fei{13k+6}$ is stably birational.
\qed
\enddemo

\proclaim{Theorem 2.3.3} Let $X$ be a nonsingular projective threefold of general type and 
suppose $P_k(X)\ge 3$, then $\fei{m}$ is birational for all $m\ge 10k+8$.
\endproclaim
\demo{Proof} When $\dim\fei{k}(X)\ge 2$, we know from Case 1 of Theorem 2.3.2 that
$\fei{m}$ is birational for $m\ge 9k+4$. When $|kK_X|$ is composed of a pencil,
from the proof of Theorem 2.3.1, we see that $\fei{k}$ will derive a fibration
$f:X'\lrw W$ onto a nonsingular curve. If $b:=g(W)>0$, then $\fei{m}$ is birational 
for $m\ge 2k+4$. 

The remain case is the one when $b=0$.  We
have an injection $\Cal{O}(2)\hookrightarrow f_*\omega_{X'}^k$. So, for each $p>0$,
we have
$$\Cal{O}(1)\otimes f_*\omega_{X'/{\Bbb P}^1}^p=
\Cal{O}(2p+1)\otimes f_*\omega_{X'}^p
\hookrightarrow f_*\omega_{X'}^{k(p+1)+p}.$$
Thus Koll\'ar's method tells that $\phi_{6k+5}$ is birational, $\phi_{4k+3}$
is generically finite and that $\dim\phi_{3k+2}(X)\ge 2$.
Now using our method, we can see that $mK_{X'}$ is effective for $m\ge 4k+4$
by Lemma 1.7. Since $(4k+3)K_{X'}$ is also effective,
thus $\phi_m$ is birational for $m\ge 10k+8$.
\qed
\enddemo

\proclaim{Corollary 2.3.2} Let $X$ be a nonsingular projective threefold of general type and
suppose $p_g(X)\ge 3$, then $\fei{m}$ is birational for $m\ge 11$.
\endproclaim
\demo{Proof} Keep the same notations as in the proof of Theorem 2.3.2.
When $\dim\fei{1}(X)\ge 2$, we set $L_3:=4K_{X'}$, $L_2=L_1:=K_{X'}$. Then 
$|L_3|$ gives a generically finite map by virtue of Case 1, Theorem 2.3.2. 
Using Lemma 1.5, we see that $|K_{X'}+2L_3+L_2+L_1|$ gives a birational map.
Thus $\fei{11}$ is birational.

When $\dim\fei{1}(X)=1$, we see from the proof of Theorem 2.3.3 that $\fei{11}$
is also birational.
\qed
\enddemo

Theorem 2.3.1, Theorem 2.3.2, Theorem 2.3.3 and Corollary 2.3.2 imply the main 
theorem.

\head {3. Open problems}\endhead
\subhead 3.1\endsubhead
Let $X$ be a nonsingular projective variety of general type of dimension $n$.
We define

$k_0(X):=min\{k|\ P_k(X)\ge 2\};$

$k_s(X):=min\{k|\ \phi_m$ is birational for $m\ge k\};$

$\mu_s(X):=\frac{k_s(X)}{k_0(X)}$, which is called {\it the relative pluricanonical
stability} of $X$. Obviously, $\mu_s(X)$ is a birational invariant.

$\mu_s(n):=\text{sup}\{\mu_s(X)|\ X$ is a $n$-fold of general type$\}$, which is
called the $n$-th {\it relative pluricanonical stability}.

It is well-known that $\mu_s(1)=3$ and $\mu_s(2)=5$ (\cite{1}). {}From the main theorem,
we have $\mu_s(3)\le 16$. What is the exact value of $\mu_s(3)$? It is also interesting to study
$\mu_s(n)$ for $n\ge 4$, even we don't know whether we should have
$\mu_s(n)<+\infty$.

\subhead 3.2\endsubhead
We would like to ask a very natural question which never happens in surface case.

\proclaim{Question} Does there exist a smooth projective threefold $X$ of general type
and two positive integers $k_1<k_2$ such that $\fei{k_1}$ is birational while
$\fei{k_2}$ is not birational?
\endproclaim

Of course, it may happen for some threefold that $P_{k_1}>P_{k_2}$ even if 
$k_1<k_2$. But we have not found any counter example yet to the above question.

\head {References}\endhead
\roster
\item"[1]" W. Barth, C. Peter, A. Van de Ven, {\it Compact Complex Surface},
Springer-Verlag 1984.
\item"[2]" A. R. Fletcher, {\it Contributions to Riemann-Roch in projective
3-folds with only canonical singularities and application},
Proc. Sympos. Pure Math. {\bf 46}, Amer. Math. Soc. Providence, 1987, 221-232.
\item"[3]" Y. Kawamata, {\it A generalization of Kodaira-Ramanujam's
vanishing theorem},  Math. Ann., {\bf 261}(1982), 43-46.
\item"[4]" Y. Kawamata, K. Matsuda, K. Matsuki, {\it Introduction to the 
minimal model problem}, Adv. Stud. Pure Math. {\bf 10} (1987), 283-360.
\item"[5]" J. Koll\'{a}r, {\it Higher direct images of dualizing
sheaves I},  Ann. of Math. {\bf 123}(1986), 11-42.
\item"[6]" J. Koll\'ar, S. Mori, {\it Birational geometry of algebraic 
varieties}, Cambridge Univ. Press, 1998.
\item"[7]" M. Reid, {\it Canonical 3-folds}, Journ\'ees de G\'eom\'etrie 
Alg\'ebrique d'Angers, A. Beauville (editor), Sijthoff and Noordhoff, Alphen aan
den Rijn, 1980, pp. 273-310.
\item"[8]" I. Reider, {\it Vector bundles of rank 2 and linear systems on algebraic
surfaces}, Ann. of Math. {\bf 127}(1988), 309-316.
\item"[9]" S. G. Tankeev, {\it On n-dimensional canonically polarized varieties and 
varieties of fundamental type}, Izv. A. N. SSSR, S\'er. Math. {\bf 35}(1971), 
31-44.
\item"[10]" E. Viehweg, {\it Vanishing theorems},  J. reine angew. Math.,
{\bf 335}(1982), 1-8.
\item"[11]" G. Xiao, {\it Finitude de l'application bicanonique des surfaces de type
g\'en\'eral}, Bull. Soc. Math. France {\bf 113}(1985), 23-51.
\endroster

\enddocument